\documentclass[reqno]{amsart}

\usepackage{amsmath,amsfonts,amssymb,latexsym}
\usepackage{hyperref}

\usepackage{times}
\usepackage{euscript} 
\DeclareMathOperator{\Ext}{Ext}

\DeclareMathOperator{\Def}{Def}
\DeclareMathOperator{\Sing}{Sing}

\def\ra{\rightarrow}
\def\cal{\mathcal} 
\def\wt{\widetilde}

\def\CC{\mathbb{C}}

\def\QQ{\mathbb{Q}}
\def\ZZ{\mathbb{Z}}

\def\OO{\cal O} 
\def\XX{\cal X}

\def\yy{\cal Y}

\def\ss{\Sigma}
\def\T{\Theta}

\def\bcp{\mathbb C\mathbb P}

\newtheorem{main}{Theorem}

\newtheorem{thm}{Theorem}[section]
\newtheorem{prop}[thm]{Proposition}

\newtheorem{defn}[thm]{Definition}

\newtheorem{lem}[thm]{Lemma}

\newtheorem{rmk}[thm]{Remark}

\numberwithin{equation}{section}

\newtheorem{ack}{Acknowledgments}        
\makeatother

\begin{document}

\title{The algebraic rational blow-down}

\author{Rare{\c s} R{\u a}sdeaconu}

\address{Rare{\c s} R{\u a}sdeaconu,
         Department of Mathematics, University of Michigan, 2074 East Hall,  530 Church Street, 
         Ann Arbor,  48104, MI, USA}

\email{rares@umich.edu}

\author{Ioana {\c S}uvaina}

\address{Ioana {\c S}uvaina, 
Department of Mathematics, 
         SUNY at Stony Brook, 
         Stony Brook, 11794, NY, USA}
        
\email{sio@math.sunysb.edu}



\date{\today}

\begin{abstract}
The normal connected sum construction of Gompf and the rational blowing-down 
technique of Fintushel - Stern are important tools in constructing symplectic 4-manifolds. 
In some cases, the 4-manifolds created this way are of K{\" a}hler type. In this 
article we investigate the occurrence of this phenomenon and give relevant examples. 
\end{abstract}

\maketitle

\bigskip

\section*{Introduction}
\label{intro}

A symplectic manifold is a smooth compact manifold endowed with a non-degenerate closed 2-form. 
For a long time, the only known rich source of examples was given by the class of K\"ahler 
manifolds. In fact, there were known only a few scattered examples of non-K\"ahler symplectic manifolds. 
The major breakthrough in constructing non-K\"ahler symplectic manifolds was realized by Gompf 
\cite{gompf}, who introduced the normal connected sum, the gluing of two symplectic manifolds along 
diffeomorphic submanifolds satisfying a suitable compatibility condition. Another source of symplectic 
manifolds was later introduced by Fintushel and Stern \cite{fs}. It is a surgery procedure, 
called the rational blow-down. It is specific to the 4-dimensional realm and amounts 
to removing a neighborhood of a linear chain of embedded spheres and replacing it with a rational ball. 
The fact that it can be performed within the class of symplectic manifolds is a result of 
M. Symington \cite{syg}, \cite{syg2}.       

It is a difficult problem to decide the existence or the non-existence of complex structures on 
symplectic 4-manifolds constructed via these two techniques. The question we would like to address here 
is the following: if the initial manifold admits an integrable complex structure, when  does the resulting
 manifold admit a integrable complex structure?
% when they actually yield complex surfaces.
 We provide a sufficient condition:
% to achieve complex structures on the rational blowing down: 
\begin{main}
\label{criterion}
Let $G$ be a finite group acting with only isolated fixed points on a smooth, compact, complex surface 
$S$ with $H^2(S,\T_S)=0.$ If the singularities of $S{/G}$ are of class $T,$ then the full rational 
blowing down $\wt S$ of the minimal resolution of $S{/G}$ admits complex structures. Moreover, 
as a smooth $4-$manifold, $\wt S$ is oriented diffeomorphic to the generic fiber of a 1-parameter $\QQ-$Gorenstein smoothing of $S{/G}.$
\end{main}

To explain our result, we should recall that the singularities of type $T$ are either 
rational double points or quotient singularities of a particular type. The exceptional 
divisor of the minimal resolution of this last type of quotient singularities is a 
linear chain of rational curves on which the rational blowing down can be performed. 
What we mean by the full rational blowing down is rationally blowing down {\em all} of 
the exceptional divisors which appear resolving the singularities which are not ordinary 
double points.

As an application we study a series of examples obtained by letting $\ZZ_4,$ the multiplicative group of roots of order $4$ of unity, act on a 
product of two appropriated Riemann surfaces. As a particular case, we find a complex structure on 
an  example constructed by Gompf, on which the existence of complex structures was still an open 
problem.  

To briefly recall Gompf's example, we start with a simply connected, relatively minimal, elliptic 
surface, with no multiple fibers, and with Euler characteristic $c_2=48.$ Up to diffeomorphisms
 \cite{elliptic}, there is only one such elliptic surface, which we call $E(4).$ It is known that
 $E(4)$ admits at most nine rational $(-4)-$curves as disjoint sections. We can form the normal
 connected sum of it with $n$ copies of $\bcp_2,$ the complex projective plane,  identifying a conic in each $\bcp_2$ with one
 $(-4)-$curve of $E(4).$ This operation 
is the same as rationally blowing-down $n ~ (-4)-$curves, and we obtain Gompf's examples 
\cite[page 564]{gompf} denoted by $W_{4,n},$ where $n=1,\dots ,9.$ The manifold $W_{4,1}$ is does not
admit any complex structure, as it violates the Noether Inequality, but the existence of complex
 structures is topologically unobstructed in the other cases.
 In Gompf's paper, he shows that the $4-$manifolds  $W_{4,n},$ for 
$n=2,3,4,$ and 9 are diffeomorphic to complex surfaces. We find a complex structure on $W_{4,8}.$ 

In this article we discuss the above examples from the point of view of deformation theory. We prove:  
\begin{main}
\label{exgompf}
The $4-$manifolds $W_{4,n}, n=2,3,4,8,9$ admit a complex structure.
\end{main}
%We should point out that the complex structure on $W_{4,8}$ is new.

The emphasis is on the methods we employ. One technique comes from Manetti's interpretation of 
the rational blowing down in algebraic setting as the 1-parameter $\QQ-$Goren-stein smoothing of a 
certain class of normal surface singularities. The second method is a natural approach towards the normal 
connected sum construction and it consists in viewing it as the smoothing of a simple normal crossing 
algebraic variety.

\section{Two constructions of symplectic 4-manifolds}

In this section we briefly recall two important construction procedures of
symplectic $4-$manifolds. One is the normal connected sum procedure of
Gompf \cite{gompf} and the other is Fintushel and Stern's \cite{fs} rational 
blow-down. In the end we describe the $4-$manifolds $W_{4,n}.$

\subsection{The normal connected sum}
Let $(M_i, \omega_i),~i=1,2$ be two symplectic $4-$manifolds. Suppose there exist
$\ss_i\subset M_i,~i=1,2$ two closed, smooth $2-$dimensional
symplectic submanifolds of the same genus, and satisfying the compatibility condition:
$$
N_{\ss_1|M_1}=N_{\ss_2|M_2}^{\vee}.
$$ 
Let $N_1(\ss_i),~N_2(\ss_i)$ be two tubular neighborhoods of $\ss_i,$ such that 
${\overline {N_1(\ss_i)}}\subset N_2(\ss_i).$ We denote by $W_i$ the tubular shell neighborhoods 
$N_2(\ss_i)-{\overline {N_1(\ss_i)}}$ of $\ss_i$ in $M_i.$ Suppose we have an orientation 
preserving diffeomorphism $\Phi:W_1\ra W_2$ taking the inside boundary of $W_1$ to the outside 
boundary of $W_2.$ We define the {\em normal connected sum} of $M_1$ and $M_2$ along $\ss_1$ 
and $\ss_2$ via $\Phi$ to be the smooth oriented manifold obtained by gluing 
$M_1-{\overline {N_1(\ss_1)}}$ and $M_2-{\overline {N_1(\ss_2)}}$ along the tubular shell 
neighborhoods $W_1$ and $W_2$ using $\Phi.$ Let $M=M_1\#_{\Phi}M_2$ be the resulting $4-$manifold. 
\begin{thm}[Gompf] Possibly after rescaling $\omega_1$ or $\omega_2,$ there 
exists a symplectomorphism $\Phi$ of the tubular shell neighborhoods of $\ss_1$ and $\ss_2$ such 
that the $4-$manifold $M=M_1\#_{\Phi} M_2$ admits a symplectic form which agrees  
with the rescaled symplectic forms on $M_i-{\overline {N_1(\ss_i)}}.$ 
\end{thm}  
%A particular case is the fiber sum operation. We start with a complex surface $S(1)$ 
%fibering over a smooth curve $C$ and choose two smooth fibers $F_1$ and $F_2.$ 
%We can form the normal connected sum of two copies of $S$ identifying $F_1$ and $F_2,$ 
%via the identity. We obtain a new $4-$manifold $S(2)\ra C \# C$ fibering over the usual connected sum of $C$ with itself, denoted $C\# C.$ 
%Of course, $S(2)$ satisfies the hypothesis of the fiber sum construction and we can iterate the procedure.  
% 
% \begin{xpl}
% \label{e(n)}
% Choose two cubics $C_1$ and $C_2$ in $\PP_2.$ The pencil generated by these two cubics has $9$ base 
% points which we blow-up. We denote by $E(1)$ the resulting rational surface. $E(1)$ will fiber over 
% $\PP_1$ with elliptic curves as fibers, and we can repeatedly apply the fiber sum operation on
% it. We obtain in this way  sequence of smooth, symplectic 4-manifolds $E(n)$ mentioned in 
% introduction. We should point out that they all admit complex structures as elliptic fibrations 
% over $\PP_1$ with no multiple fibers.   
% \end{xpl}

\subsection{The rational blow-down}
Let $C_{p,q}$ be a open smooth $4$-manifold obtained by plumbing disk bundles
over the configuration of $2$-sphere prescribed by the linear graph:

\begin{picture}(400,40)(-100,-20)
  \put(0,3){\makebox(200,20)[bl]{\hspace{7pt}$b_{k}$ \hspace{21pt} $b_{k-1}$
                                  \ \ \hspace{89pt}  $b_{1}$}}
\multiput(10,0)(40,0){2}{\line(1,0){40}}
                         \multiput(10,0)(40,0){3}{\circle*{3}}
                         \multiput(100,0)(5,0){4}{\makebox(0,0){$\cdots$}}
  \put(125,0){\circle*{3}}
  \put(125,0){\line(1,0){40}}
  \put(165,0){\circle*{3}}
\end{picture}

\noindent
Here $p>q>0$ are two relatively prime integers and  $\frac{p^{2}}{pq-1}
=[b_{k},b_{k-1}, \ldots, b_{1}]$ is the unique continued fraction with all
$b_{i} \geq 2.$ Each vertex represents a disk bundle over the 
$2$-sphere of self-intersection $-b_{i}.$ Then $C_{p,q}$ is a negative definite simply connected 
$4-$manifold whose boundary is a lens space $L(p^2, 1-pq).$ The lens space 
$L(p^2, 1-pq)$ bounds a rational ball $B_{p,q}$ with $\pi_{1}(B_{p,q}) \cong {\ZZ}_{p}.$ 

%Furthermore, the inclusion
%$\partial B_{p,q} \longrightarrow B_{p,q}$ induces an epimorphism
%$$
%\pi_{1}(\partial B_{p,q})\cong {\mathbf Z}_{p^2} \longrightarrow \pi_{1}(B_{p,q}) \cong {\mathbf Z}_{p}.
%$$ 
%
% 
Suppose $X$ is a smooth $4$-manifold containing a configuration $C_{p,q}$. Then we may construct 
a new smooth $4$-manifold $X_{p,q}$, called the {\em (generalized) rational blow-down} of $X.$ by 
replacing $C_{p,q}$ with the rational ball $B_{p,q}.~X_{p,q}$ is uniquely determined (up to 
diffeomorphism) from $X,$ as  each diffeomorphism of $\partial B_{p,q}$ extends over the rational ball 
$B_{p,q}.$ It is proved in \cite{syg2} that if $C_{p,q}$ is symplectically embedded in $X,$ then the 
rational blow-down carries a symplectic structure.

\section{An algebraic description of the rational blow-down}

The point of view we adopt in this paper is that the generalized rational blowing down procedure 
and the smoothing of isolated complex surface singularities are essentially the same in the algebraic setting. While in this article we merely give some immediate applications of this idea, we will explore it in more depth in a forthcoming article. Here we consider only the case of some quotient singularities, a case well documented in the literature.
 
Following Manneti's presentation \cite{man} of the results of \cite{ksb}, we begin by recalling the terminology and some general results.  
\begin{defn}
A normal variety $X$ is $\QQ-$Gorenstein if it is Cohen-Macaulay and a multiple of the canonical 
divisor is Cartier.
\end{defn}
\begin{defn}
A flat map $\pi:\XX\ra \CC$ is called a one-parameter $\QQ-$Gorenstein smoothing of a
 normal 
singularity $(X,x)$ if $\pi^{-1}(0)=X$ and there exists $U\subset \CC$ an open neighborhood of $0$ 
such that the following conditions are satisfied. 
\begin{itemize}
\item [ i)] $\XX$ is $\QQ-$Gorenstein
\item [ ii)] The induced map $\XX\ra U$ is surjective 
\item [ iii)] $X_t=\pi^{-1}(t)$ is smooth for every $t\in U-\{0\}.$
\end{itemize}
\end{defn}
\begin{defn}
A normal surface singularity is \emph {of class $T$} if it is a quotient singularity and admits a 
$\QQ-$Gorenstein one-parameter smoothing. 
\end{defn}
The following result of Koll{\'a}r and Shepherd-Barron \cite {ksb} gives a complete description of 
the singularities of class $T.$
\begin{prop}[Koll{\'a}r,~Shepherd-Barron]
\label{descsing}
The singularities of class $T$ are the following:
\begin{itemize}
% \item [ 1)] Smooth points
\item [ 1)] Rational double points;
\item [ 2)] Cyclic singularities of type $\displaystyle \frac{1}{dn^2}(1,dna-1),$ for $d>0,~n\geq 2$ 
and $(a,n)=1.$ 
\end{itemize}
\end{prop}

Let $a,d,n>0$ be integers with $(a,n)=1,$ and $\yy\subset \CC^3\times \CC^d$ be the 
hypersurface of equation 
$$
uv-y^{dn}=\sum_{k=0}^{d-1}t_ky^{kn},
$$ 
where $t_0,\dots ,t_{d-1}$ are linear coordinates on $\CC^d.~\ZZ_n$ acts on $\yy,$ the 
action being generated by:
$$
(u,v,y,t_0,\dots ,t_{d-1})\mapsto (\zeta u,\zeta ^{-1}v,\zeta ^{a}y,t_0,\dots ,t_{d-1})
$$
Let $\XX=\yy{/\ZZ_n}$ and $\phi:\XX\ra \CC^d$ the quotient of the projection $\yy \ra \CC^d.$
\begin{prop}
\label{defT} 
$\phi:\XX\ra \CC^d$ is a $\QQ-$Gorenstein deformation of the cyclic singularity $(X,x)$ of type 
$\displaystyle \frac{1}{dn^2}(1,dna-1).$

Moreover, every $\QQ-$Gorenstein deformation $\XX\ra \CC$ of a singularity $(X,x)$ of type 
$\displaystyle \frac{1}{dn^2}(1,dna-1)$ is isomorphic to the pullback of $\phi$ for some germ of 
holomorphic map $(\CC,0)\ra(\CC ^d,0).$ 
\end{prop}
The following proposition due to Manetti \cite{man} gives an algebraic interpretation of the algebraic 
rational blow-down.
\begin{prop}
\label{arbd}
Let $X$ be a compact complex surface with singularities of class $T,~\wt X$ be its minimal resolution 
and $\pi:\XX\ra U$ a $\QQ-$Gorenstein smoothing of $X.$ Then the full rational blow-down of $\wt X$ 
is oriented diffeomorphic to $X_t=\pi^{-1}(t),$ for any $~t\neq 0.$   
\end{prop}

Manetti's algebraic description of the rational-blowing-down procedure gets more substance by 
providing a global criterion for smoothings of singularities. 
% 
% This is however well-know, e.g  
% \cite{wahl} or the more recent \cite{mancrit}.    

\smallskip

Let $X$ be a compact, reduced analytic space with $\Sing (X)=\{x_1,\dots , x_n\}.$ By restriction, for 
each $i=1,\dots , n,$ any deformation of $X$ defines a deformation of the singularity $(X,x_i)$ 
with \emph {the same base space.} Thus, if $\Def X$ and $\Def (X,x_i)$ denote the base space of the 
versal deformation of $X$ and $(X,x_i),$ respectively, there exists a natural morphism: 
$$
\Phi:\Def X\ra \prod_{i=1}^{n} \Def (X,x_i)
$$

If all of the singularities $(X,x_i)$ are of class $T,$ for each $i$ we can choose in  $\Def (X,x_i)$ the (smooth) component 
corresponding the $\QQ-$Gorenstein deformations.
What we are interested in here is when these deformations of  singularities lift to a deformation 
of the total space. The answer is given by a general criterion of Wahl: 
\begin{prop}
\label{wman}
Let $\T_X$ be the tangent sheaf of $X.$ If $H^2(X,\T_X)=0,$ then the morphism $\Phi$ is smooth. 
In particular, every deformation of the singularities $(X,x_i),~i=i,\dots ,n$ may be globalized. 
\end{prop}
For the proof of this proposition we refer the interested reader to either \cite[page 242]{wahl} or 
to \cite[page 93]{mancrit}.

This criterion becomes especially useful in the following situation. 

\smallskip

Suppose we start with a smooth compact complex surface $X.$ Let $G$ be a finite group acting on 
$X$ with isolated fixed points. Let $Y=X/G$ be the quotient, $\{y_1,\dots ,y_n\}=\Sing(Y)$ the 
singular locus of $Y,$ and $f:X\ra Y$ the quotient map. 
We denote by $\T_X$ the holomorphic tangent bundle of $X,$ and by $\T_Y={(\Omega_Y^1)}^\vee,$ the tangent sheaf of $Y.$ 

\begin{lem} 
\label{quotcrit}
In the above notations, $H^2(X,\T_X)=0\Longrightarrow H^2(Y,\T_Y)=0.$
\end{lem}

\begin{proof} To prove the vanishing of $H^2(Y,\T_Y)=0,$ we need to have a convenient description 
of the tangent sheaf $\T_Y$ of $Y.$ In our case, this is provided by Schlessinger \cite{schl}: 
$$
\T_Y=(f_*\T_X)^G.
$$
Now, by averaging, we get a map $f_*\T_X\ra\T_Y=(f_*\T_X)^G.$ But this means $\T_Y$ is a direct summand 
of $f_*\T_X.$ To finish the proof, since $f$ is a finite map, the Leray spectral sequence provides 
an isomorphism:
$$
H^2(X,\T_X)\simeq H^2(Y,f_*\T_X),
$$
and the conclusion of the lemma follows.
\end{proof}

\begin{proof}[Proof of Theorem \ref{criterion}] Assume now $G$ acts on a smooth complex surface $S$ with fixed points only. If the singularities 
of $S{/G}$ are of class $T$ only, we can look at the components of each versal deformation space of 
any such singular point, and pick the one corresponding to $\QQ-$Gorenstein deformations.  
Theorem \ref{criterion} follows immediately from the algebraic description of the rational blowing down, 
the globalization criterion \ref{wman} and Lemma \ref{quotcrit}.
\end{proof}

\begin{rmk}
We should point out that the smoothings of rational double points are diffeomorphic to their minimal
resolution. Thus, for the singularities of class $T$ the full rational 
blowing down, essentially performed only on the minimal resolution of
the quotient singularities which are not rational double points, coincides with the simultaneous smoothing of all the singular points. 
\end{rmk}

\section{A family of examples}

In this section, by adopting the above viewpoint, we will exhibit a complex structure on Gompf's example $W_{4,8}.$ 

\smallskip

Let $C\subset \bcp_1\times \bcp_1$ be the smooth curve of genus $3$ given by the equation 
$$
F({\bf z}, {\bf w})=z_0^4(w_0^2+w_1^2)+z_1^4(w_0^2-w_1^2)=0, 
$$ 
where $({\bf z}, {\bf w})=([z_0:z_1],[w_0:w_1])$ are the standard bi-homogeneous coordinates 
on $\bcp_1\times \bcp_1.$ 

An easy calculation shows that the action of the cyclic group $\ZZ_4$ on $\bcp_1\times \bcp_1$ 
generated by 
$$
([z_0:z_1],[w_0:w_1])\mapsto ([iz_0:z_1],[w_0:w_1])
$$
has four fixed points of $C:$ two points $P_1=([0:1],[1:1])$ and $P_2=([0:1],[1:-1]),$ where in local 
coordinates $\ZZ_4$ acts by multiplication with $i,$ and two $Q_1=([1:0],[1:i])$ and  
$Q_2=([1:0],[1:-i]),$ where in local coordinates our group acts by multiplication with $-i.$ 

% The group $\ZZ_4 \oplus \ZZ_4$ acts on $ C \times C$. 
We will be interested in the manifold obtained by taking the quotient under the diagonal action of $\ZZ_4$. Let $X= (C \times C)/{\ZZ_4}.$
% Consider now the singular complex surface:
% $$
% X=(C\times C)/{\ZZ_4}\hookrightarrow (\PP_1\times \PP_1 \times \PP_1 \times \PP_1)/{\ZZ_4}.
% $$
This action has $16$ fixed points. At eight of them, $(P_i,P_j),$ and 
$(Q_i,Q_j),~i=1,2$ the group $\ZZ_4$ acts (in local coordinates) as 
$$
(z_1,z_2)\mapsto (iz_1,iz_2),
$$ 
while at the other eight $(P_i,Q_j)$ and $(Q_i,P_j),~i=1,2$ it acts as   
$$
(z_1,z_2)\mapsto (iz_1,-iz_2).
$$ 
Thus, the singular complex 2-dimensional variety $X$ will have 8 singular points of type $A_3$ 
and 8 quotient singularities of type $\frac{1}{4}(1,1).$ The minimal resolution of the last type 
of singularities consists in replacing each such singular point by a smooth rational curve of self-intersection $(-4).$ Let 
$\widehat X$ be the minimal resolution of $X.$
\begin{prop}
\label{gex}
$\widehat X$ is a simply connected, minimal, elliptic complex surface with no multiple fibers and
with the Euler characteristic $c_2=48.$
\end{prop}
\begin{proof}
The quotient $C/ \ZZ_4$ is a rational curve, and we will denote by $P'_1, P'_2, Q'_1$ and $Q'_2$ the image of $P_1, P_2, Q_1$ and $Q_2,$ respectively, under the projection map. 
Let: 
$$
\pi_1 : \widehat X\longrightarrow  C/ \ZZ_4\cong \bcp_1
$$
the factorization of the projection on the first factor. This fibration has four
singular fibers above $P'_i$ and $Q'_i,~i=1,2.$ The generic fiber is
a Riemann surface of genus $3$, while each of the singular
fibers consist of a chain of nine spheres, one of which is
the quotient $C/ \ZZ_4$ and the other eight are exceptional
spheres introduced by the resolution of singularities. It follows 
that our manifold is simply connected.

To see the elliptic fibration, we consider first the covering:
$$
 \pi_1 \times \pi_2 : \widehat X \longrightarrow (C \times
 C)/ (\ZZ_4 \oplus \ZZ_4) = C/ \ZZ_4 \times C/\ZZ_4 \cong
 \bcp_1\times \bcp_1.$$
As the construction is symmetric in the two factors, we will
 use freely the identification of $C/ \ZZ_4$ with
 $\bcp_1$ with four marked points $P'_i,Q'_i,~i=1,2.$ Let 
$B,~D\subset\bcp_1\times \bcp_1$ be the following divisors: 
\begin{align*}
B=&~\pi_1^{-1}(P'_1)\cup
\pi_1^{-1}(P'_2)\cup\pi_2^{-1}(Q'_1)\cup\pi_2^{-1}(Q'_2)\\
D=&~\pi_1^{-1}(Q'_1)\cup
 \pi_1^{-1}(Q'_2)\cup\pi_2^{-1}(P'_1)\cup\pi_2^{-1}(P'_2).
\end{align*} 
$\widehat X$ can be seen as a
bi-double cover of $\bcp_1\times \bcp_1,$ first branched over the union of $B$ and
 $D$, then branched over the union of the total transforms of $B$
 and $D.$ 
 
Remark that $\OO_{\bcp_1\times\bcp_1}(B) \cong\OO_{\bcp_1\times\bcp_1}(D)=\OO_{\bcp_1\times\bcp_1}(2,2)$ and the
 generic element of the associated linear system is a smooth elliptic curve.
 Let $L$ be the pencil generated by $B$ and $D.$ The base
 locus of this pencil, $B \cap D =\lbrace(P'_i,P'_j),
 (Q'_i,Q'_j),~i,j=1,2\rbrace$, is a subset of  the set of singular points
 of the branch locus of the
 first double cover. The linear system $(\pi_1 \times \pi_2)^* (L)$ will be base
 point free, with smooth elliptic curve as general members. This gives us
 the elliptic fibration, containing no rational curves of self-intersection $(-1).$ 
 The exceptional divisors introduced
 above $B\cap D$ are the eight sections of the elliptic fibration, all of self-intersection $-4.$
 
% The same $4:1$ map $\pi_1 \times \pi_2 : \widehat X \ra \bcp_1\times \bcp_1,$ can be used to immediately check the minimality of $\widehat X.$
Next, we need to know two topological invariants, the Euler characteristic $\chi(\widehat X)$ and signature $\sigma(\widehat X),$ for example. Both computations are immediate if we look at the quotient map $C\times C \longrightarrow X$. We obtain $\chi(\widehat X)=48$ 
and $\sigma (\widehat X)= -32$, which imply $c_1^2(\widehat X)=0.$ This also yields the minimality of $\widehat X.$

\end{proof}

In particular, by \cite{elliptic}, $\widehat X$ is diffeomorphic to $E(4).$ 
Applying Theorem \ref{criterion}, we can conclude now the existence of complex 
structures on $W_{4,8}.$ 

\bigskip

Our example can be easily generalized in the following way. Consider the smooth curves
 $C_k,~C_l\subset\bcp_1\times \bcp_1$ given by the equations
$$
z_0^4f_k(w_0,w_1)+z_1^4g_k(w_0,w_1)=0
$$  
and 
$$
z_0^4f_l(w_0,w_1)+z_1^4g_l(w_0,w_1)=0, 
$$  
respectively. Here $(f_k,~g_k)$ and $(f_l,~g_l)$ are generic pairs of homogeneous polynomials of
 degree $k,$ and $l,$  respectively. The above discussion can be now easily repeated for $C_k\times
 C_l$ and the induced $\ZZ_4$ action.

\section{The algebraic normal connected sum}

In this section we approach the normal connected sum procedure from the algebraic point of view. 
We will test our point of view on Gompf's examples $W_{4,n},$ for $n=2,3,4$ and $9.$ 

\smallskip

Mimicking the symplectic normal sum, we start with two pairs  of complex varieties 
$(X_1, Y_1),~(X_2,Y_2),$ where $X_i,~i=1,2$ are smooth and $Y_i\subset X_i,~i=1,2$ are smooth 
subvarieties satisfying the following conditions:
\begin{itemize}
\item $Y_1\simeq Y_2;$
\item $N_{Y_1|X_1}=N_{Y_2|X_2}^{\vee}.$
\end{itemize}
Using, the condition on the normal bundles, we can glue $X_1$ and $X_2$ to form a normal crossing complex variety $X.$ Then the symplectic normal 
sum can be interpreted as a smoothing of $X,$ in the sense of Friedman \cite{friedman}, as long as the 
smoothing is of K{\" a}hler type. Before we proceed, we recall some basic facts on the deformation 
theory of singular spaces.    

\smallskip 

\subsection{Deformation theory of normal crossing varieties} 
\label{gen}

Let $X$ be a smooth complex manifold, and $Y\subset X$ a smooth submanifold. 
By $\T_X$ we denote the sheaf of holomorphic vector fields of $X,$ and $\T_{X,Y}$ will 
denote the sheaf of holomorphic vector fields on $X$ which are tangent to $Y.$ 

If $Z$ is a compact, singular, reduced complex space, the deformation theory of normal crossing 
varieties \cite{friedman} is given in terms of the global $\Ext$ groups 
$T_Z^i:=\Ext^i(\Omega_Z, \OO_Z).$ $T_Z^1$ will describe the infinitesimal deformations of the complex 
structure of $Z,$ and the obstructions lie in $T_Z^2.$ These groups are usually computed from 
their "local" versions, the sheaves ${\cal E}xt^i(\Omega_Z, \OO_Z)=: \tau_Z^i,$ using the "local to 
global" spectral sequence $E_2^{p,q}=H^p(\tau^q_Z)\Rightarrow T_Z^{p+q}.$ 

Suppose $Z$ is a simple normal crossing variety, i.e. $Z$ locally looks like a union of hyperplanes 
and whose irreducible components are smooth. In this case, the local to global spectral sequence gives:
\begin{equation}
\label{l2g}
0\ra H^1(\tau^0_Z)\ra T_Z^1\ra H^0(\tau^1_Z)\ra H^2(\tau^0_Z)\ra T_Z^2\ra 0.
\end{equation}
The space $H^1(\tau^0_Z)$ classifies all "locally trivial" deformations of $X,$ i.e. for which the 
singularities remain locally a product. Their obstructions lie in $H^2(\tau^0_Z).$  

In the cases treated below, $Z$ is chosen to be $d-$semistable, that is  $T_Z^1=\OO_D,$ where 
$D=\Sing(X)$ is the singular locus of $X.$  

\begin{defn}
We say that a proper flat map $\pi:\XX\ra \Delta$ from a smooth $(n+1)-$fold $\XX$ to 
$\Delta=\{|z|<1\}\subset \CC$ is a smoothing of a reduced, not necessarily irreducible complex 
analytic variety $X,$ if $\pi ^{-1}(0)=X$ and $\pi ^{-1}(t)$ is smooth for $t\in \Delta$ sufficiently 
small.
\end{defn}
\begin{rmk} 
\label{crit}
\emph {Before we proceed with our examples, we should point out that 
if $H^2(\tau ^0_Z)= H^1(\OO_D)=0,$ then $T_Z^2=0$, and so the deformation 
problem is unobstructed. If $X$ is $d-$ semistable, it will admit \cite{friedman} a versal 
(one-parameter) deformation with smooth total space, and smooth generic fiber. In this case, 
we say we have a one-parameter smoothing of $X.$}
\end{rmk}

\subsection{Gompf's examples} 
In what follows we are going to treat a very particular situation, and give a simple cohomolgical criterion suitable to the study of Gompf's examples. 

\medskip

Let $S$ be a complex surface, containing $n$ smooth, disjoint, rational curves of self-intersection 
$-4,$ which we are going to denote by $D_1,\dots ,D_n.$ As discussed, for each of these 
curves, we can glue in a copy of $\bcp_2$ to form a simple normal crossing complex 
variety denoted by $Z$ with $n+1$ irreducible components, and 
$\Sing (Z)=D_1+\cdots +D_n.$ Since each copy of $\bcp_2$ is glued along a conic, it 
follows that $Z$ satisfies the $d-$semistability condition. An easy, but useful criterion is:
\begin{prop} 
\label{sncs}
$Z$ admits a one-parameter smoothing if 
$$
H^2(S,\T_S\otimes \OO_S(-D_1-\cdots-D_n))=0.
$$ 
\end{prop}
\begin{proof}
From Remark \ref{crit}, it suffices to check whether $H^2(\tau ^0_Z)= H^1(D,\OO_D)=0,$ where $D=D_1+\cdots +D_n.$ 

Since the $D_i$'s are smooth, disjoint, rational curves, it follows that $H^1(D,\OO_D)=0.$ The sheaf $\tau ^0_Z,$ naturally sits in the 
exact sequence: 
\begin{equation}
\label{se}
0\ra \T_{S}\otimes \OO_S(-D)\oplus \bigoplus_{i=1}^n \T_{\bcp_2}\otimes\OO_{\bcp^2}(-C_i)\ra \tau ^0_Z \ra \T_{D}\to 0,
\end{equation}
where by $C_i$ we denoted the smooth conic in $\bcp_2$ corresponding to $D_i.$ But, for any smooth conic $C\subset \bcp_2$ we have $H^2(\bcp_2,\T_{\bcp_2}\otimes \OO_{\bcp_2}(-C))=0.$ A simple inspection of the cohomology sequence associated to 
(\ref{se}) concludes the proof.
\end{proof}

We take now a look from our perspective at Gompf's examples $W_{4,n},$ for $n=2,3,4$ or 
$9.$ In these cases, he found complex structures as appropriate multiple covers of 
$\bcp_2.$ 
Guided by his complex structures, what we do here is merely to reprove this in an algebraic,  more conceptual way, 
using the method described above. We will discuss only the $W_{4,2}$ case, the rest of the cases 
follow analogously, and are left to the reader.  

\smallskip

We start with the Hirzebruch surface $\ss_4.$ Let $f:X\ra \ss_4$ be the double cover of $\ss_4,$ 
branched along a smooth member of the linear system $D\in|4(C_0+4f)|.$ Such a smooth member exists, 
as a consequence of the standard results on linear systems on Hirzebruch 
surfaces \cite{hart}. Here $C_0$ is the negative section and $f$ is the class of a fiber. 
$X$ is a smooth, simply connected elliptic surface, diffeomorphic to $E(4).$ Moreover, since $D$ and 
$C_0$ are disjoint, $X$ contains exactly two smooth rational curves of self-intersection 
$-4,$ the two irreducible components of the preimage of $C_0.$ We denoted these two curves 
by $D_1$ and $D_2.$  

We perform now the algebraic normal connected sum along these two curves gluing in two copies of 
$\bcp_2,$ each along a smooth conic, and denote the newly formed singular variety by $Z$. 
From Proposition \ref{sncs} we know that if $H^2(X,\T_X\otimes \OO_X(-D_1-D_2))=0,$  there is no obstruction to the smoothing of $Z$ . But, using the structure of $X$ as a 
double covering of $\ss_4,$ and the Leray spectral sequence, we get:
\begin{align*}
&~H^2(X,\T_X\otimes \OO_X(-D_1-D_2)) \\
=&~H^2(X,\T_X\otimes f^*\OO_{\ss_4}(-C_0))\\
=&~H^2(\ss_4,f_*\T_X\otimes \OO_{\ss_4}(-C_0))\\
=&~H^2(\ss_4,\T_{\ss_4}\otimes \OO_{\ss_4}(-C_0))\oplus H^2(\ss_4,\T_{\ss_4}\otimes \OO_{\ss_4}(-3C_0-8f)).
\end{align*}
To prove the vanishing of the last two cohomology groups, we use the Serre duality: 
\begin{align*}
&~H^2(\ss_4,\T_{\ss_4}\otimes \OO_{\ss_4}(-C_0))\\
=&~H^0(\ss_4,\Omega^1_{\ss_4}\otimes \OO_{\ss_4}(C_0)\otimes 
\OO_{\ss_4}(K_{\ss_4}))\\
=&~H^0(\ss_4,\Omega^1_{\ss_4}\otimes \OO_{\ss_4}(-C_0-6f)).
\end{align*}
Here $K_{\ss_4}=-2C_0-6f$ denotes the canonical divisor of $\ss_4.$ Since the divisor $C_0+6f$ is effective and is linearly equivalent to a smooth curve \cite{hart}, we have the exact sequence of sheaves: 
$$
0\ra \Omega^1_{\ss_4}\otimes \OO_{\ss_4}(-C_0-6f))\ra \Omega_{\ss_4}^1 \ra 
\OO_{C_0+6f}\ra 0.
$$
Passing to the cohomology sequence and since $H^0(\ss_4,\Omega_{\ss_4}^1)=0,$ 
we get the vanishing of $H^2(\ss_4,\T_{\ss_4}\otimes \OO_{\ss_4}(-C_0)).$ 

For the vanishing of the second term we proceed in the same fashion:
\begin{align*}
&~H^2(\ss_4,\T_{\ss_4}\otimes \OO_{\ss_4}(-3C_0-8f))\\
=&~H^0(\ss_4,\Omega^1_{\ss_4}\otimes \OO_{\ss_4}(3C_0+8f)\otimes 
\OO_{\ss_4}(K_{\ss_4}))\\
=&~H^0(\ss_4,\Omega^1_{\ss_4}\otimes \OO_{\ss_4}(C_0+2f)).
\end{align*}
Arguing by contradiction, if there exists a global non-zero section of $\Omega^1_{\ss_4}\otimes \OO_{\ss_4}(C_0+2f),$ then it must exist global non-zero section of 
$$ \displaystyle
\bigwedge ^2 ( \Omega^1_{\ss_4}\otimes \OO_{\ss_4}(C_0+2f)) \cong \OO_{\ss_4}(K_{\ss_4})\otimes \OO_{\ss_4}(2C_0+4f)
\cong \OO_{\ss_4}(-2f).
$$
But this is impossible.

\begin{rmk}
{\emph {We believe that working a little bit harder, we should be able to prove that the complex we 
found is actually the same as the one Gompf described for $W_{4,2}$, namely as the double covering of 
$\bcp_2$ along a smooth octic, a Horikawa surface. To give some strong indications, we can prove that 
the normal connected sum from our point of view of $\bcp_2$ with $\ss_4$ glued along a 
conic and the negative section, respectively is isomorphic to $\bcp_2.$ We have a natural  
a holomorphic map $F:Z\ra \bcp_2\#_C \ss_4,$ extending the double covering of $\ss_4.$ What we found 
above are the smoothings of $Z$ and $\bcp_2\#_C \ss_4.$ Using an appropriate deformation theory, we 
should be able to show that the map $F$ deforms too, to yield the expected double covering $W_{4,2}\ra 
\bcp_2,$ branched along a smooth octic. }}
\end{rmk}

The computations for the other examples go along the same lines, with minor modifications. 
For $W_{4,3}$ we end up with the complex structure of the $3:1$ covering of $\bcp_2$ branched along a 
smooth sextic curve, $W_{4,4}$ is the simple bi-double cover of $\bcp_2$ branched along a transverse pair 
of conics, and $W_{4,9}$ is a $\ZZ_3\oplus \ZZ_3$ cover of $\bcp_2$ branched along 3 transverse conics.
With the results obtained in the previous section, the proof of Theorem \ref{exgompf} is now complete.

\begin{ack} {\emph{We thank Claude LeBrun for suggesting this problem to us, and Dusa McDuff for her useful comments. The first author would also like to thank for the hospitality of IMADA, University of Southern Denmark where most of this paper was  written.}} 
\end{ack}

\bibliographystyle{alpha}

\vskip1cm

\end{document}